\documentclass[12pt,letterpaper]{article}
\usepackage{geometry}
\usepackage[utf8]{inputenc}

\usepackage{amsmath}
\usepackage{amssymb}
\usepackage{amsthm}

\usepackage{hyperref}

\usepackage{graphicx}
\usepackage{float}
\usepackage{enumitem}
\usepackage{caption}
\usepackage{subcaption}

\usepackage{authblk}

\newcommand{\norm}[1]{\left\lvert#1\right\rvert}

\newcommand{\R}{{\mathbb{R}}}
\newcommand{\bx}{\mathbf{x}}
\newcommand{\bv}{\mathbf{v}}

\theoremstyle{definition}
\newtheorem{definition}{Definition}[section]
\theoremstyle{remark}
\newtheorem{remark}{Remark}[section]
\title{Concurrent Emergence of Clustering, Flocking and Synchronization in Systems of Interacting Agents}
\author[1]{Trenton Gerew}
\author[1]{Ming Zhong}
\affil[1]{\small{\emph{Department of Applied Mathematics, Illinois Institute of Technology}}}

\date{\today}

\begin{document}
\maketitle
\begin{abstract}
We present some recent development in modeling concurrent emergence of collective behaviors, namely, the emergence of clustering, flocking and synchronization at the same time.  We derive two new models, namely Swarmalator-Vicsek and Swarmalator-Cucker-Smale, which can produce the synchronization and flocking with interesting spatial patterns.  We present extensive numerical insights into how the synchronization of phases can affect both the spatial patterns and the flocking behavior.
\end{abstract}
\section{Introduction}
Collective behaviors (also known as self organization), such as clustering \cite{Krause2009, BHT2009, MT2014}, flocking \cite{starling, murmurations, birds, flap}, milling \cite{PhysRevLett.96.104302}, swarming \cite{PhysRevLett.120.198101}, synchronization\cite{kuramoto, winfree, strogatz, sync}, describes how global orders, can emerge from initially randomized configuration through local interactions between pairs of agents.  It has gathered considerable amount of research interests lately, see the Section $9$ in \cite{MT2014} and references therein. Moreover it also has a great number of of applications in Physics (such as self assemble nanoparticles, superconductivity, crystal growth), Chemistry (molecular self-assembly, self-assembled monolayers), Biology (morphogenesis, social behavior of insects, chemotaxis/phototaxis), human society (herd behavior, self-referentiality, self-organizing market economy), etc.  With proper modeling tools, single kind of behavior can be replicated through carefully designed interacting agent systems.  We present in this paper that a couple new models that can reproduce complicated behaviors with the mixture of three collective behaviors, namely synchronization, flocking and clustering, which is partially motivated by the synchronized fireflies behavior observed in the nature \cite{sync, geometry, swarmalator}.

One of the first models which combines two collective behaviors together is the swarmalator model proposed in \cite{swarmalator}.  This swarmalator model combines clustering and synchronization together.  To be precise, it is about the how synchronization of phases affects the spatial pattern of the agents.  The system is described by an ODE system, i.e., for a system of $N$ interacting agents, their states are evolved according to the following ODE equations,
\[
\begin{cases}
\dot\bx_i    &= \sum_{j = 1, j \neq i}^N\phi^E(\norm{\bx_j - \bx_i}, \xi_j - \xi_i)(\bx_j - \bx_i), \\
\dot\xi_i &= \sum_{j = 1, j \neq i}^N\phi^{\xi}(\norm{\bx_j - \bx_i}, \xi_j - \xi_i)(\xi_j - \xi_i),
\end{cases} \quad i = 1, \cdots, N.
\]
Here the pair $(\bx_i, \xi_i) \in \R^d\times\R$ presents the position/phase for the $i^{th}$ agent respectively in the systems.  Furthermore, the functions, $\phi^E, \phi^{\xi}:\R^+\times\R \rightarrow \R$, are called the interaction functions, which govern how the $j^{th}$ agent influences the change of position and phase of the $i^{th}$ agent.  Moreover, they depend only on the pairwise distance between the two agents as well as the relative phase difference.  Although whole system has local (pairwise) interactions between the agents, it can produce interesting spatial patterns with various phase distribution among the agents (static or non-static with either synchronized distribution of phases), see \cite{swarmalator} for detailed examples.  

However, individuals in nature tend to align and move as others do in their neighborhood \cite{OBRIEN19891}.  The synchronized flapping of migratory birds in the V-shape \cite{birds}, \cite{flap}, synchronized firefights (most significant sighting in the smokey mountain area) are intriguing examples of concurrent emergence of more than two collective behaviors.  Hence the combination of only two behaviors is not enough to model natural phenomena.  Therefore we propose two new models, where the the states of the agents now include $(\bx_i, \xi_i, \theta_i)$ or $(\bx_i, \xi_i, \bv_i)$, where $\theta_i$ is the heading (the polar angle of velocity of the $i^{th}$ agent) and $\bv_i$ is the velocity of the $i^{th}$ agent respectively.  The ODE system is changed to incorporate the additional state.  For example, the Swarmalator-Viscek model is governed by the following ODE system
\[
\begin{cases}
\dot\bx_i &= \bv_i + \begin{bmatrix}\cos(\theta_i) \\ \sin(\theta_i)\end{bmatrix} + \frac{1}{\norm{\Lambda_i}}\sum_{j \in \Lambda_i} \phi^E(\norm{\bx_j - \bx_i}, \theta_j - \theta_i, \xi_j - \xi_i)(\bx_j - \bx_i), \\
\dot\theta_i &= \frac{1}{\norm{\Lambda_i}}\sum_{j \in \Lambda_i}\phi^{\theta}(\norm{\bx_j - \bx_i}, \theta_j - \theta_i, \xi_j - \xi_i)(\theta_j - \theta_i), \\
\dot\xi_i &= \omega_i + \frac{1}{\norm{\Lambda_i}}\sum_{j \in \Lambda_i} \phi^{\xi}(\norm{\bx_j - \bx_i}, \theta_j - \theta_i, \xi_j - \xi_i)(\xi_j - \xi_i),
\end{cases} \quad i = 1, \cdots, N.
\]
This new system will have the states, namely position/heading angle/phase, for the $i^{th}$ agent, denoted as $(\bx_i, \theta_i, \xi_i)$.  Furthermore, the interaction of these internal states together with states from their neighboring agents can produce complex behaviors, especially the concurrent emergence of clustering, flocking and synchronization.    With the inclusion of three different states for each individual agent, these two new models can now combine the flocking behavior described in \cite{vicsek, cs-model} to synchronization of phases and spatial patterns described in \cite{swarmalator}.  We are able to show through extensive numerical experiments that when phases also enter into the interaction with position and velocity, intriguing spatial patterns can emerge.  For detailed discussion as well as the model equations for the Swarmalator-Cucker-Smale model, see section \ref{sec:model}.

The remaining of this paper is arranged as follows: Section \ref{sec:model} describes in details how the two new models are derived; Section \ref{sec:numerics} shows detailed experiments of the different patterns which our new models can produce; and we conclude our paper in Section \ref{sec:conclude}.
\section{Model Description}\label{sec:model}
We describe in this section in details about the two new models, namely Swarmalator-Viscek and Swarmalator-Cucker-Smale.  The two models govern how the state tuples, i.e. $(\bx_i, \xi_i, \theta_i)$ or $(\bx_i, \xi_i, \bv_i)$, would interact with each other and with those from neighboring agents, and such intertwined relationship for these state variables would lead to interesting flocking and synchronization with complicated spatial patterns.
\subsection{Swarmalator-Viscek Model}
First, we consider a second-order swarmalator model with headings, i.e. $\theta_i \in \mathcal{S}^2$ (the $2D$ unit sphere), for each agent.  The governing equations for a population of $N$ swarmalators are
\begin{equation}\label{eqn:swarmalator_viscek}
\begin{cases}
\dot\bx_i &= \bv_0 + \begin{bmatrix}\cos(\theta_i) \\ \sin(\theta_i)\end{bmatrix} + \frac{1}{\norm{\Lambda_i}}\sum_{j \in \Lambda_i} \phi^E(\norm{\bx_j - \bx_i}, \theta_j - \theta_i, \xi_j - \xi_i)(\bx_j - \bx_i), \\
\dot\theta_i &= \frac{1}{\norm{\Lambda_i}}\sum_{j \in \Lambda_i}\phi^{\theta}(\norm{\bx_j - \bx_i}, \theta_j - \theta_i, \xi_j - \xi_i)(\theta_j - \theta_i), \\
\dot\xi_i &= \omega_0 + \frac{1}{\norm{\Lambda_i}}\sum_{j \in \Lambda_i} \phi^{\xi}(\norm{\bx_j - \bx_i}, \theta_j - \theta_i, \xi_j - \xi_i)(\xi_j - \xi_i),
\end{cases} \quad i = 1, \cdots, N.
\end{equation}
Here the tuple $(\bx_i, \xi_i, \theta_i)$ describes the position/phase/heading (in terms of polar angle of the velocity vector) of the $i^{th}$ agent respectively.  Moreover, $\bv_0, \omega_0$ are constants describing its self-propulsion velocity and natural frequency, respectively.  The functions $(\phi^E, \phi^{\xi}, \phi^{\theta})$ define the interaction between agent $j$ and agent $i$.  $\Lambda_i$ is the neighborhood around agent $i$, and for the purposes of this paper, we define the neighborhood of interaction as in \cite{vicsek}:
\begin{definition}[Neighborhood of Interaction]
Given an interaction radius $r > 0$, the neighborhood of interaction for swarmalator $i$ is
\begin{equation}\label{neighborhood}
\Lambda_i = \left\{j \in [N] \mid \norm{\bx_j - \bx_i} \leq r \right\}, \quad [N] = \{1, \cdots, N\}.
\end{equation}
\end{definition}
\begin{remark}
Note that the choice of the neighborhood parameter $r > 0$ will affect the spatial patterns significantly, see section \ref{sec:numerics} for details.
\end{remark}
When compared to the original second-order swarmalator model \cite{OBRIEN19891}, our mode has the local dependence on neighborhood to both the position and phase equations.  It is a more natural modeling of the swarming oscillators, since the decision of neighborhood for each agent will have effect on all three states instead of only on the local heading.  In fact, a more precision description of the swarming oscillators should include three different kinds of neighborhood for each state equation.  However, we will then introduce too many changing parameters in this model.  We will leave such discussion in the future work.  Moreover, the dependence of the interaction functions on pairwise state variables can produce more complex patterns.  Difference choices of $(\phi^E, \phi^{\xi}, \phi^{\theta})$ can also lead to interesting patterns, see \cite{swarmalator}.  Moreover, these interaction functions can be also data-driven i.e. learned from observation data \cite{LZTM2019, MMQZ2021, FMMZ2022}.
\subsection{Swarmalator-Cucker-Smale Model}
An important step in modeling self organization for local models is in choosing a proper neighborhood.  For the Swarmalator-Viscek model in particular, it is crucial on how to choose the radius of interaction, i.e. $r > 0$ in \eqref{eqn:swarmalator_viscek}.  As demonstrated in section \ref{sec:numerics}, the choice of $r$ has a significant influence on the behavior of the whole system.  We address the issue by switching the emphasis of local neighborhood from selection of a subset of agents to proper choices of interaction functions, i.e. the interaction function on the alignment of velocities of agents (similar to that in \cite{cs-model}).  Therefore, we present the following Swarmalator-Cucker-Smale model, for $i = 1, \cdots, N$,
\begin{equation}\label{eqn:swarmalator_cucker_smale}
\begin{cases}
\dot\bx_i &= \bv_i, \\
\dot\bv_i &= \frac{1}{N}\sum_{j = 1, j \neq i}^N\Big[\phi^E(\norm{\bx_j - \bx_i}, \norm{\bv_j - \bv_i}, \xi_j - \xi_i)(\bx_j - \bx_i) \\
&\quad + \phi^{\bv}(\norm{\bx_j - \bx_i}, \norm{\bv_j - \bv_i}, \xi_j - \xi_i)(\bv_j - \bv_i)\Big] \\
\dot\xi_i &= \omega_i + \frac{1}{N}\sum_{j = 1, j \neq i}^N\phi^{\xi}(\norm{\bx_j - \bx_i}, \norm{\bv_j - \bv_i}, \xi_j - \xi_i)(\xi_j - \xi_i), \\
\end{cases}
\end{equation}
Here $\bv_i$ is the velocity for the $i^{th}$ agent, which can change in both direction and magnitude (as compared to the Swarmalator-Viscek model, where the velocity has fixed magnitude).  The functions $(\phi^E, \phi^{\xi}, \phi^{\theta})$ can be compactly supported to reflect the local neighborhood.  

(Add more)
\section{Numerical Results}\label{sec:numerics}
We conduct numerical tests of the two models we propose in the previous sections.
\subsection{Swarmalator-Vicsek}
\begin{figure}[H]
\centering
\begin{subfigure}{.32\textwidth}
\centering
\includegraphics[width=\linewidth]{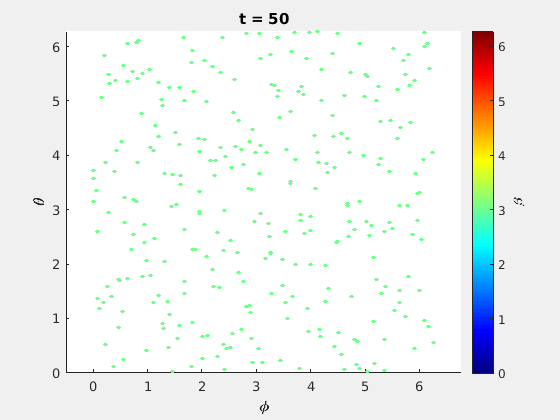}
\caption{Static asynchronous state \\ $(r,J,K) = (1.25,0.1,1)$}
\label{fig-async-scat}
\end{subfigure}%
\begin{subfigure}{.32\textwidth}
\centering
\includegraphics[width=\linewidth]{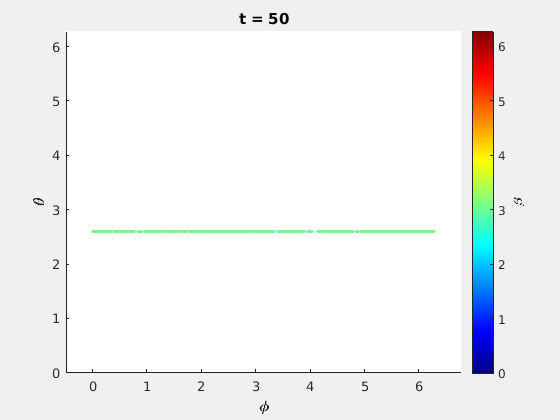}
\caption{Static gradient state \\$(r,J,K) = (1.10,0.1,-1)$}
\label{fig-grad-scat}
\end{subfigure}
\begin{subfigure}{.32\textwidth}
\centering
\includegraphics[width=\linewidth]{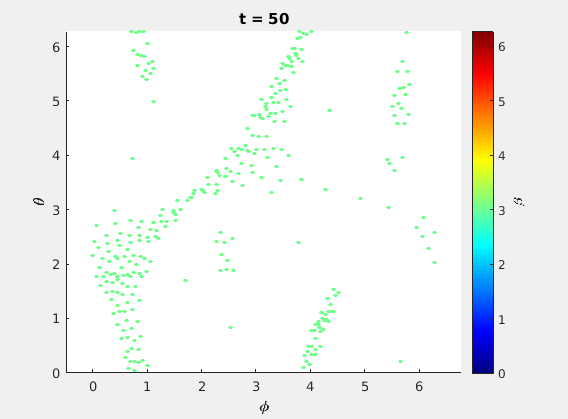}
\caption{Static wave state \\ $(r,J,K) = (1.25,1,0)$}
\label{fig-wave-scat}
\end{subfigure}%
\caption{Scatter plots in $(\psi, \xi)$ space for $N=300$ swarmalators and $T=50$ time steps, where $\psi = \tan^{-1}(y/x)$. The color of each point represents the orientation $\theta$.}
\label{fig-local01}
\end{figure}

In order to compare our results fairly to the results in \cite{swarmalator}, we choose the pair, $(\phi^E, \phi^{\xi})$, to be the following
\[
\phi^E(r, s) = \frac{A + J\cos(\xi_j - \xi_i)}{r} - \frac{B}{r^2} \quad \text{and} \quad \phi^{\xi}(r, s) = \frac{K\sin(s)}{r}.
\]
Here $r_{ij} = \norm{\bx_j - \bx_i}$ and $s_{ij} = \xi_j - \xi_i$ represents the pairwise distance and pairwise difference of phases variables.  The interaction function on the heading variable $\theta$, i.e. $\phi^{\theta}$, is taken as
\[
\phi^{\theta}(z) = 1, \quad \text{for all $z_{ij} = \theta_j - \theta_i$}.
\]
This is a continuation of the setup in \cite{vicsek}.  One can use other interaction functions.  We then set the natural frequency $\omega_i = 0$ for each agent.  Following the setup in \cite{swarmalator}, we set $A = B = 1$.  Therefore, our model has a set of three parameters, i.e. $(r, J, K)$, to choose for different types of system behaviors.  The phase coupling strength is represented by the parameter $K$.  For $K > 0$, the phase coupling tends to minimize the phase difference between swarmalators, and for $K < 0$ the phase difference is maximized.  When $K = 0$, the phases of the swarmalators do not change.  The strength at which phase similarity enhances spatial attraction is measured by parameter $J$.  If $J > 0$, swarmalators in the same phase are attracted to each other.  For $J < 0$, swarmalators move towards others with the opposite phase.  When $J = 0$, the spatial attraction is independent of phase.  We constrain $J \in [-1, 1]$ to keep $\frac{A + J\cos(\xi_j - \xi_i)}{r} - \frac{B}{r^2} \geq 0$. We performed numerical experiments to identify the behavior of our model.  All the simulations were run using MATLAB's standard ODE solver \verb|ode45|.  Initially, the swarmalators were positioned in a box of length $1$ with phases in $[-\pi, \pi]$ and orientations in $[0, 2 \pi]$, all drawn uniformly at random.  We choose $\bv_i = 0.003$ (the self propelling velocity) to be within the range of velocities explored in \cite{vicsek}.  We found the system's behavior is highly dependent on the chosen radius of interaction, and new types of grouping behavior emerge.  To compare our model to the swarmalator model, we use the parameters that describe the five unique states found in \cite{swarmalator}.
\begin{enumerate}[font=\bfseries]
  \item \textbf{Static asynchronous state.}  In the original swarmalator model the parameters $K > 0$ and all $J$ form a static synchronous phase in which the swarmalators form a circularly symmetric distribution and are fully synchronized in phase and orientation.  In our model, this combination of $J$ and $K$ form what is essentially the static asynchronous state from \cite{swarmalator}, as shown in Figure \ref{fig-async}.  The swarmalators separate faster than they can synchronize, and end up forming a crystal-like distribution in which all phases $\theta$ can occur.  At small radii $r$, no extremely little interaction occurs at all and the orientations $\beta$ are random as well.  However, after about $r = 0.65$ the interactions are strong enough that the orientations converge rapidly.  Note that in Figure \ref{fig-async-scat} the swarmalators are distributed uniformly, so that every phase occurs everywhere.
\begin{figure}[H]
\centering
\begin{subfigure}{.32\textwidth}
\centering
\includegraphics[width=\linewidth]{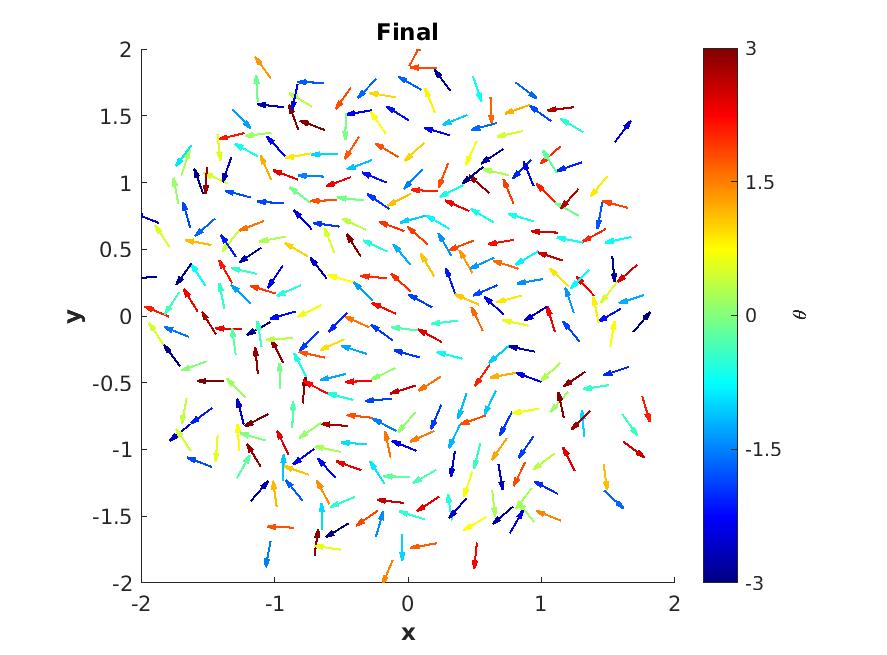}
\caption{$r = 0.20$}
\end{subfigure}%
\begin{subfigure}{.32\textwidth}
\centering
\includegraphics[width=\linewidth]{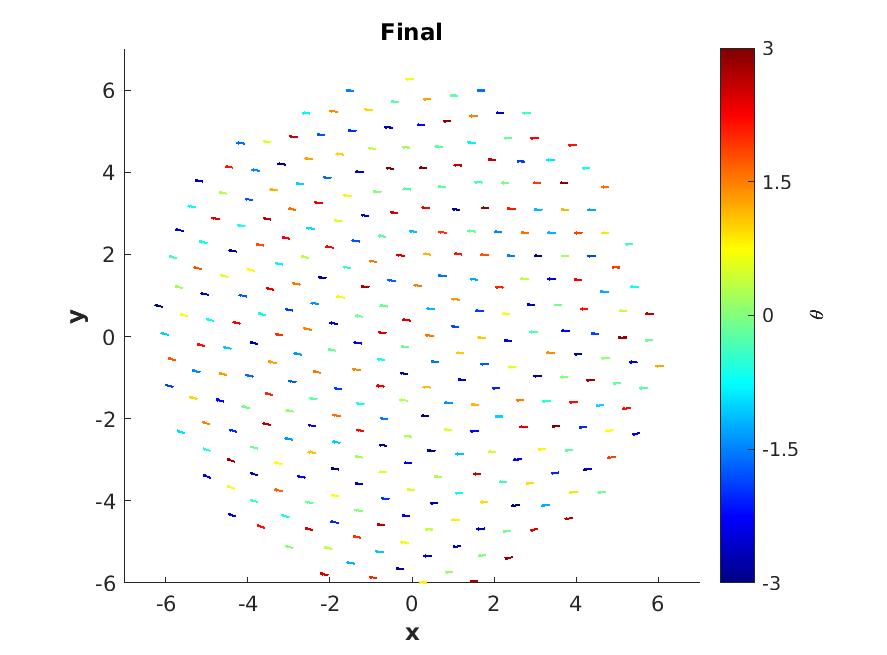}
\caption{$r = 0.65$}
\end{subfigure}
\begin{subfigure}{.32\textwidth}
\centering
\includegraphics[width=\linewidth]{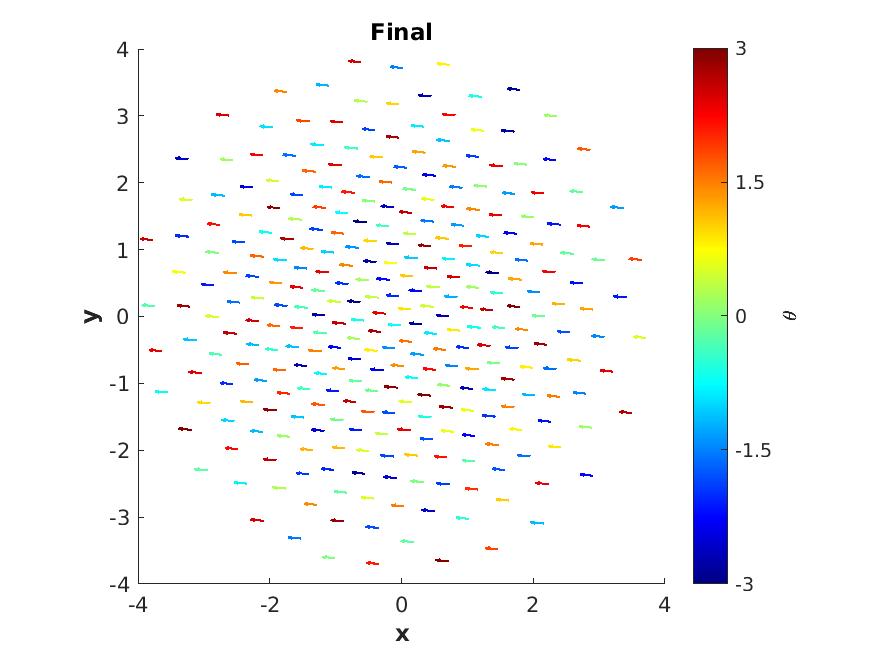}
\caption{$r = 1.40$}
\end{subfigure}
\caption{Static asynchronous state: $J = 0.1, \ K = 1$ with $N = 300$ swarmalators after $T = 50$ time units.}
\label{fig-async}
\end{figure}
\item \textbf{Static gradient state.}  Swarmalators can also form a gradient state, as seen in Figure \ref{fig-grad}.  This state occurs for $J < 0, \ K < 0$ and also for $J > 0$ if $J < \norm{K_c}$ as shown in \cite{swarmalator}. In the original swarmalator model, this parameter set forms the completely asynchronous state.
\begin{figure}[H]
\centering
\begin{subfigure}{.32\textwidth}
\centering
\includegraphics[width=\linewidth]{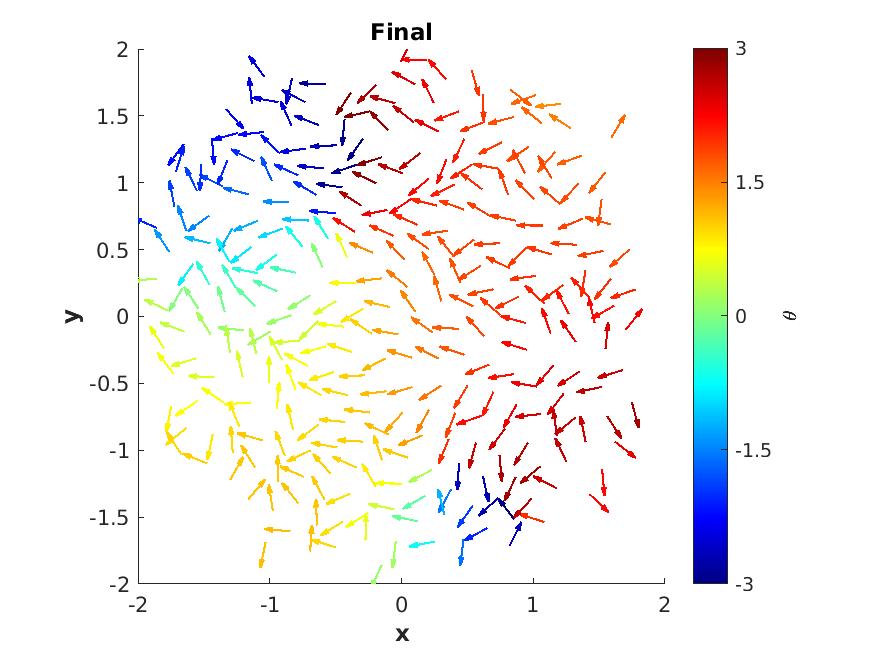}
\caption{$r = 0.20$}
\end{subfigure}%
\begin{subfigure}{.32\textwidth}
\centering
\includegraphics[width=\linewidth]{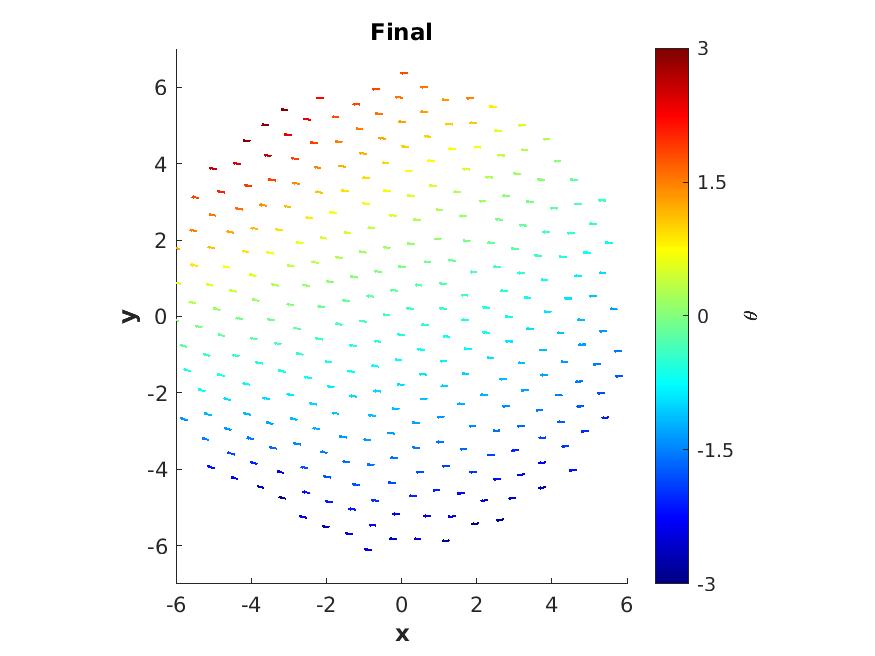}
\caption{$r = 0.65$}
\end{subfigure}
\begin{subfigure}{.32\textwidth}
\centering
\includegraphics[width=\linewidth]{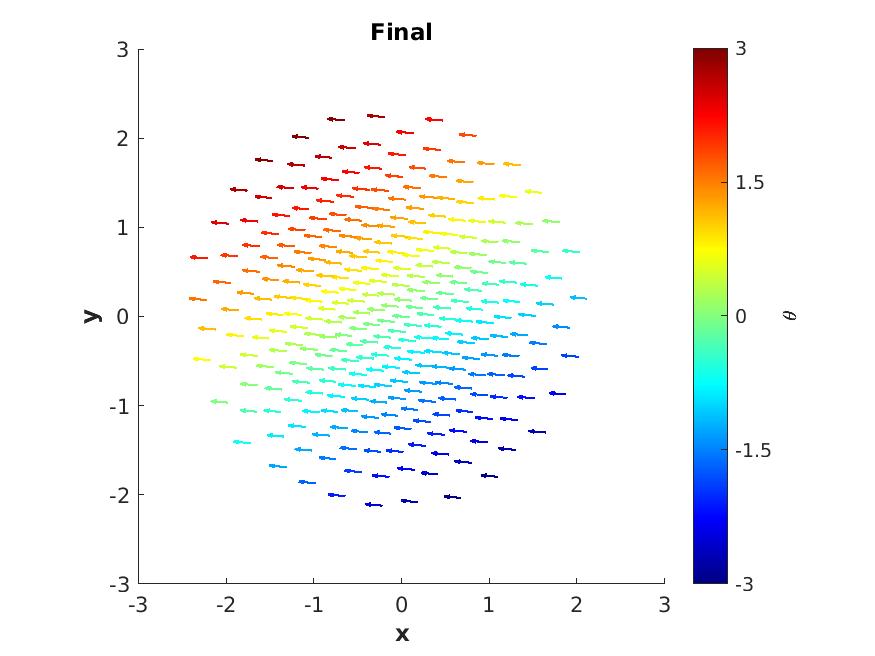}
\caption{$r = 1.40$}
\end{subfigure}
\caption{Static gradient state: $J = 0.1, \ K = -1$ with $N = 300$ swarmalators after $T = 50$ time units.}
\label{fig-grad}
\end{figure}
\item \textbf{Static wave state.}  The special case with $K = 0$ and $J > 0$ forms a curvy wave-like structure.  The swarmalators are frozen in their initial phases, and $J > 0$ causes the swarmalators to be attracted towards others with similar phases.  In the infinite-radius case, this interaction forms an annular structure called the static phase wave.  With the limited radius, however, swarmalators can't ``see'' all others of the same phase and thus settle into a repeating phase-wave-like pattern.  The new structure is not annular, but the spatial angle $\phi$ and the phase $\theta$ of each swarmalator are still moderately correlated, as seen in Figure \ref{fig-wave} and \ref{fig-wave-scat}.
\begin{figure}[H]
\centering
\begin{subfigure}{.32\textwidth}
\centering
\includegraphics[width=\linewidth]{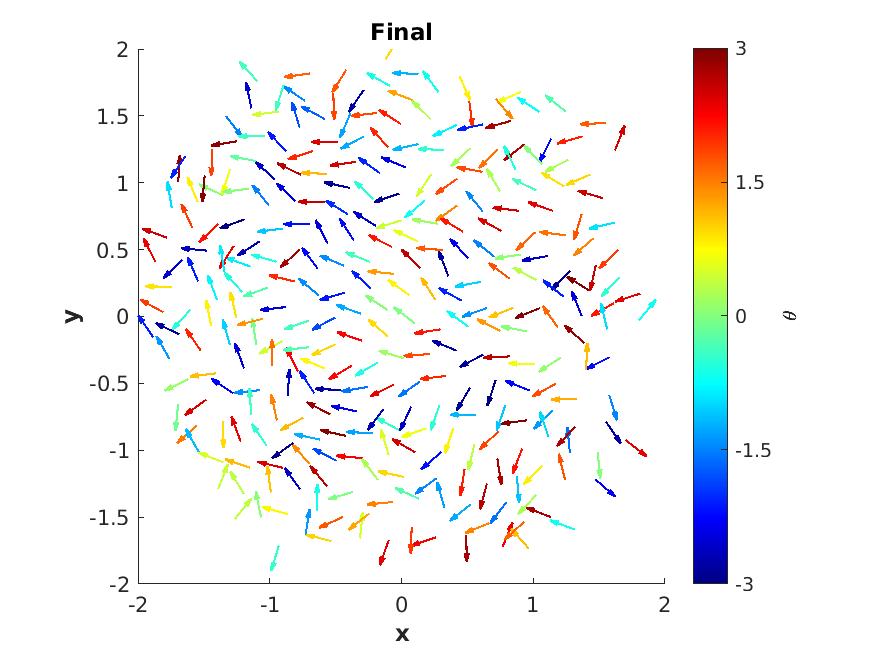}
\caption{$r = 0.20$}
\end{subfigure}%
\begin{subfigure}{.32\textwidth}
\centering
\includegraphics[width=\linewidth]{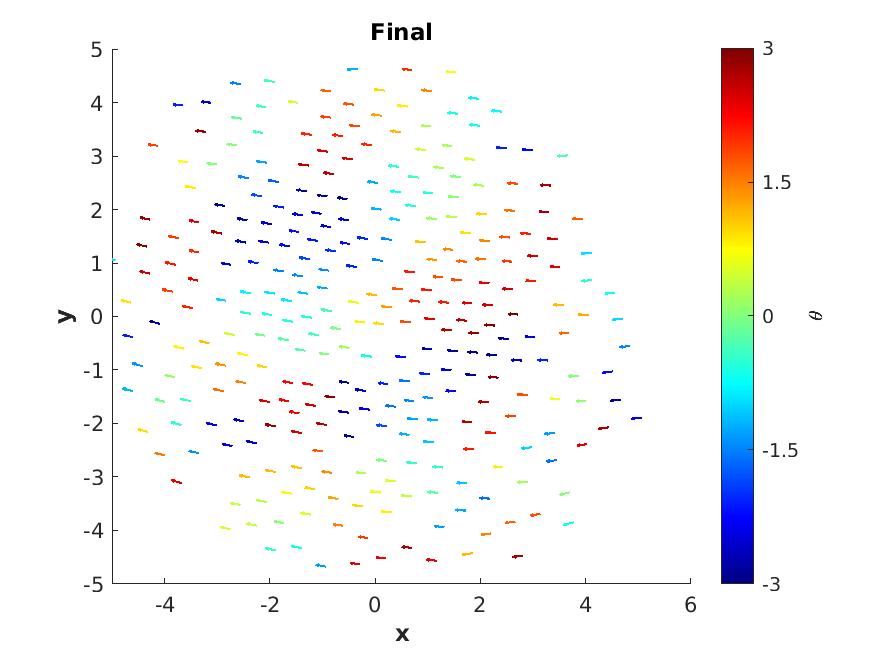}
\caption{$r = 0.65$}
\end{subfigure}
\begin{subfigure}{.32\textwidth}
\centering
\includegraphics[width=\linewidth]{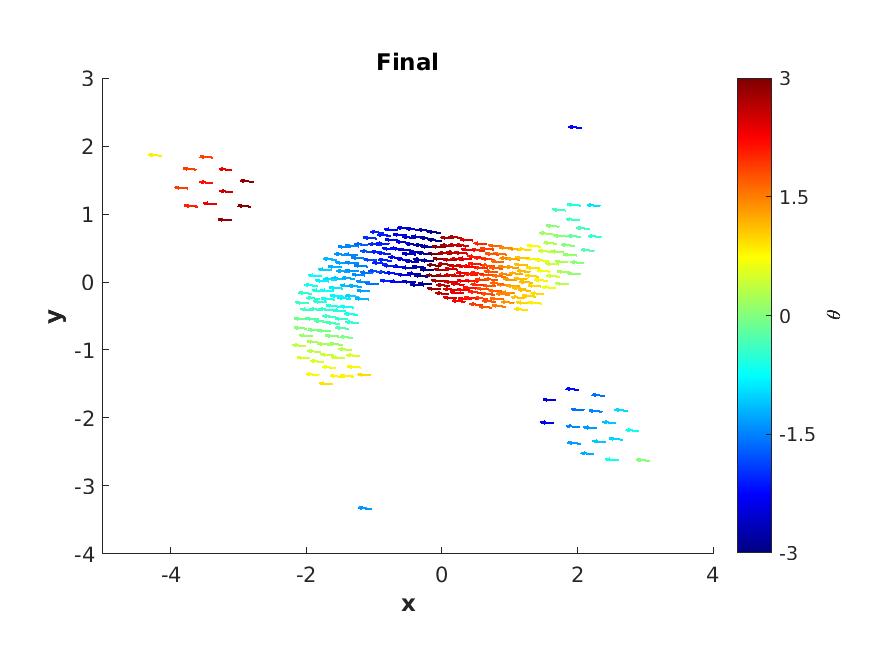}
\caption{$r = 1.40$}
\end{subfigure}
\caption{Static wave state: $J = 1, \ K = 0$ with $N = 300$ swarmalators after $T = 50$ time units.}
\label{fig-wave}
\end{figure}
\item \textbf{Clustered state.}  Moving into $K < 0$, we encounter the clustered state.  As seen in Figures \ref{fig-clusters} and \ref{fig-clust-scat}, the population splinters into groups of distinct phases.  In the movies included in the supplemental materials, it can be seen that some nearby groups will collide and merge into one phase. It is unclear what determines the number of clusters, but in general fewer are found when a larger radius is used.  This state corresponds with the splintered phase wave state in \cite{swarmalator}.
\begin{figure}[H]
\centering
\begin{subfigure}{.32\textwidth}
\centering
\includegraphics[width=\linewidth]{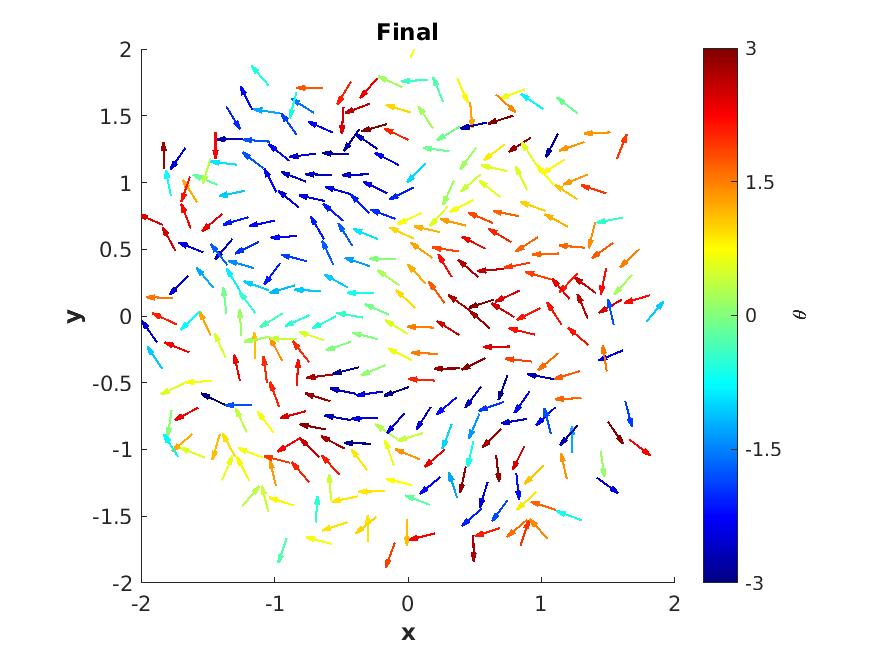}
\caption{$r = 0.20$}
\end{subfigure}%
\begin{subfigure}{.32\textwidth}
\centering
\includegraphics[width=\linewidth]{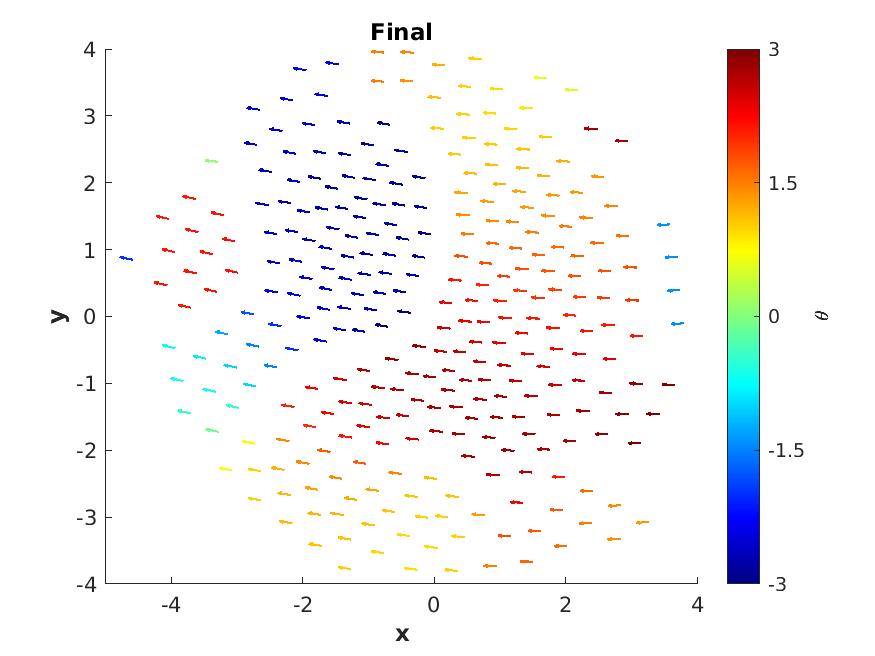}
\caption{$r = 0.65$}
\end{subfigure}
\begin{subfigure}{.32\textwidth}
\centering
\includegraphics[width=\linewidth]{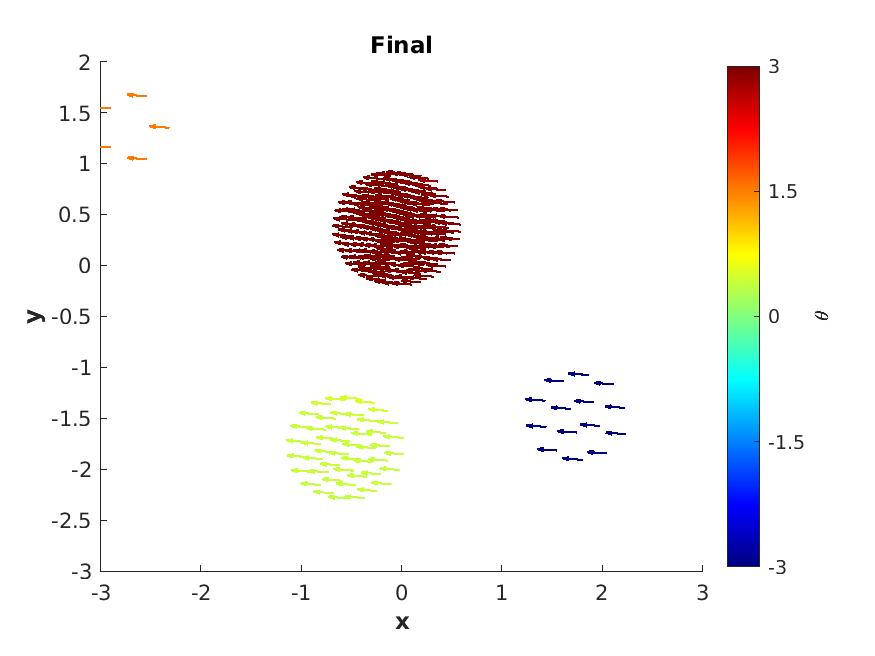}
\caption{$r = 1.40$}
\end{subfigure}
\caption{Clustered state: $J = 1, \ K = -0.1$ with $N = 300$ swarmalators after $T = 50$ time units.}
\label{fig-clusters}
\end{figure}
\item \textbf{Static synchronous state.}  The final state is the synchronous state which occurs after $K$ is further increased.  The parameter set for this state corresponds to the active phase wave for the unaltered swarmalator model.  In contrast, here the swarmalators rapidly adopt the same phase and condense into a tight group as seen in Figures \ref{fig-sync} and \ref{fig-sync-scat}.
\begin{figure}[H]
\centering
\begin{subfigure}{.32\textwidth}
\centering
\includegraphics[width=\linewidth]{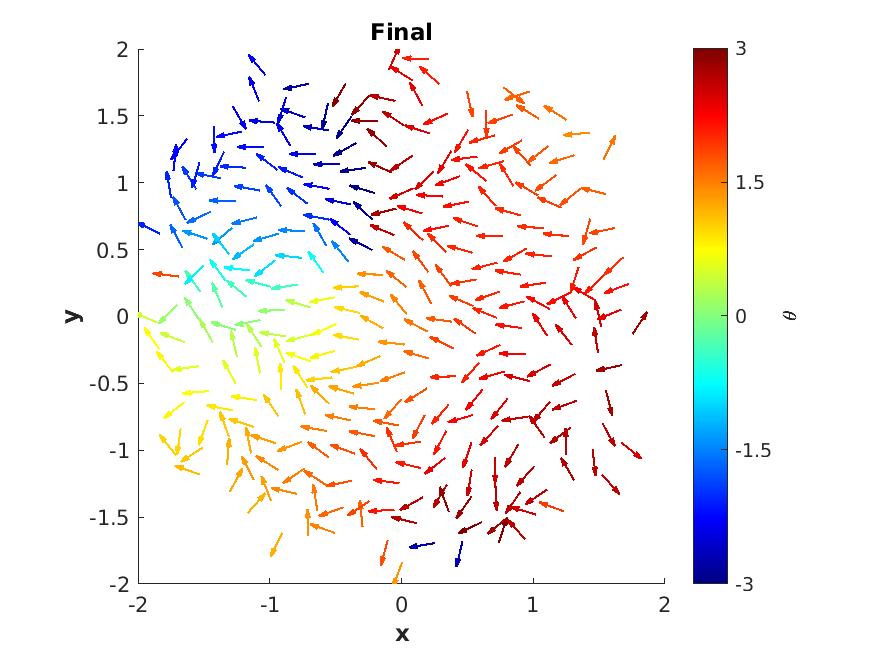}
\caption{$r = 0.20$}
\end{subfigure}%
\begin{subfigure}{.32\textwidth}
\centering
\includegraphics[width=\linewidth]{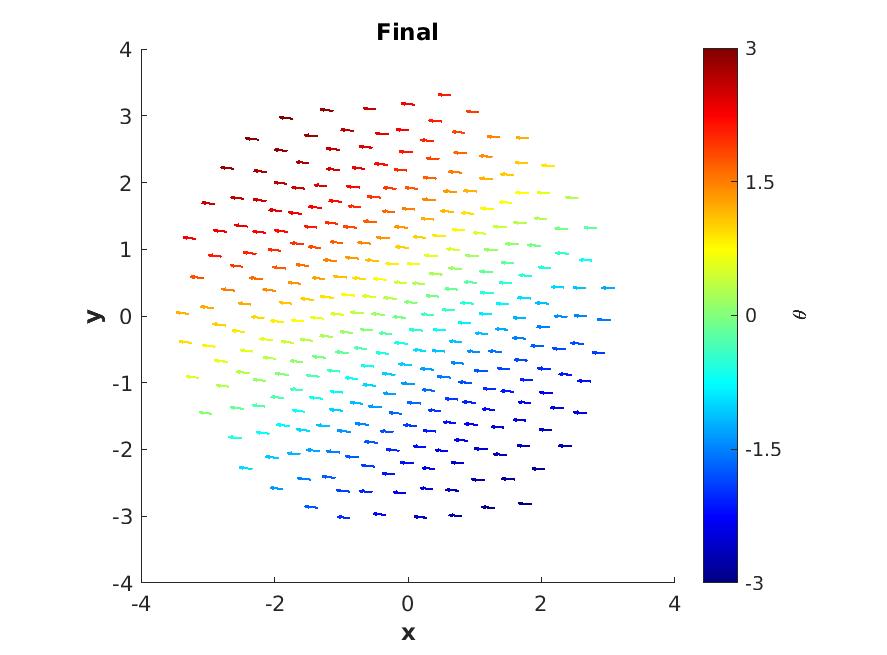}
\caption{$r = 0.65$}
\end{subfigure}
\begin{subfigure}{.32\textwidth}
\centering
\includegraphics[width=\linewidth]{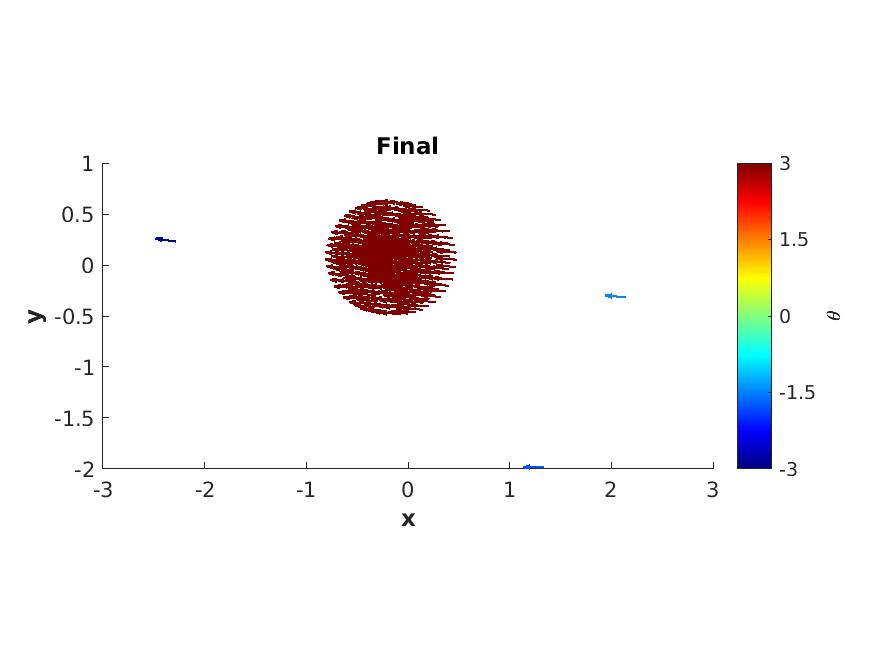}
\caption{$r = 1.40$}
\end{subfigure}
\caption{Static synchronous state: $J = 1, \ K = -0.75$ with $N = 300$ swarmalators after $T = 50$ time units.}
\label{fig-sync}
\end{figure}
\end{enumerate}
\textbf{Conclusion}:  in all cases, pattern formation does not begin to occur until about $r = 0.6$.  The reason for this is not clear, and the critical radius appears to occur for all starting densities of swarmalators.  Figure \ref{fig-local01} shows the static spatial patterns with different phases distribution for the relationship between the distribution of the polar angles of position vs the phase variable.
 Figure \ref{fig-local02} shows the astatic spatial patterns with different phases distribution for the relationship between the distribution of the polar angles of position vs the phase variable.
\begin{figure}[H]
\begin{subfigure}{.32\textwidth}
\centering
\includegraphics[width=\linewidth]{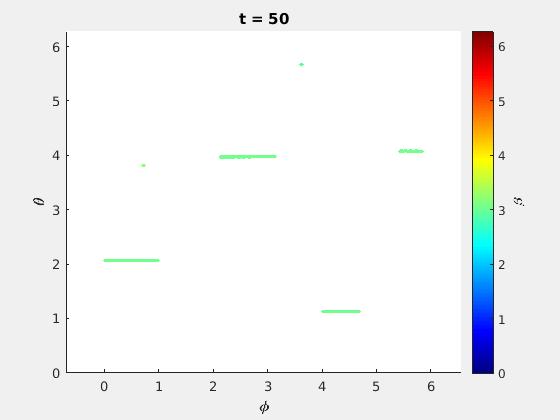}
\caption{Clustered state \\ $(r,J,K) = (1.10,1,-0.1)$}
\label{fig-clust-scat}
\end{subfigure} \hfill
\begin{subfigure}{.32\textwidth}
\centering
\includegraphics[width=\linewidth]{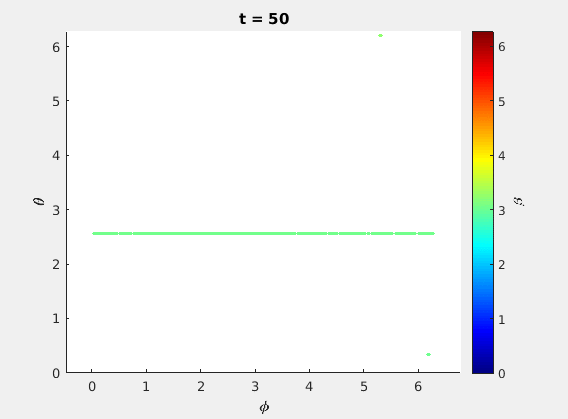}
\caption{Static synchronous state \\ $(r,J,K) = (1.25,1,-0.75)$}
\label{fig-sync-scat}
\end{subfigure}
\caption{Scatter plots in $(\psi, \xi)$ space for $N=300$ swarmalators and $T=50$ time steps, where $\psi = \tan^{-1}(y/x)$. The color of each point represents the orientation $\theta$.}
\label{fig-local02}
\end{figure}
\subsection{Swarmalator-Cucker-Smale}
The model in \eqref{eqn:swarmalator_cucker_smale} behaves essentially the same as a global swarmalator-Vicsek model posed in \cite{swarmalator}, where the parameter controlling the size of the neighborhood, i.e. $r$, is going to infinity.  We performed simulations of the model with $N = 500$ swarmalators.  The swarmalators were initially positioned in a box of length 2 with phases in $[-\pi, \pi]$ and velocity vectors in $[-0.1,0.1] \times [-0.1,0.1]$, all drawn uniformly at random.
The same three stationary states and two active states from \cite{swarmalator} are clearly identifiable.  First, the three static states are present.
\begin{figure}[H]
\centering
\begin{subfigure}{0.32\textwidth}
\centering
\includegraphics[width=\linewidth]{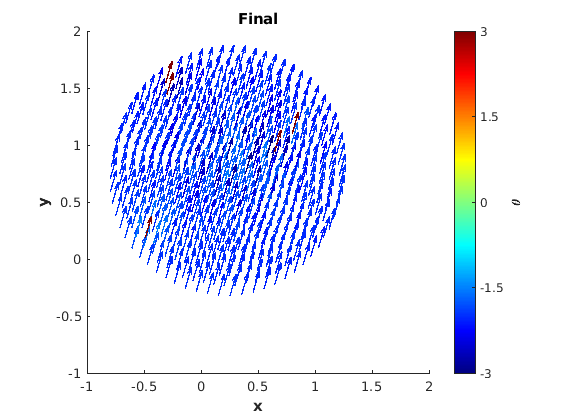}
\caption{Static sync state \\ $(J,K) = (0.1, 1)$}
\end{subfigure}%
\begin{subfigure}{0.32\textwidth}
\centering
\includegraphics[width=\linewidth]{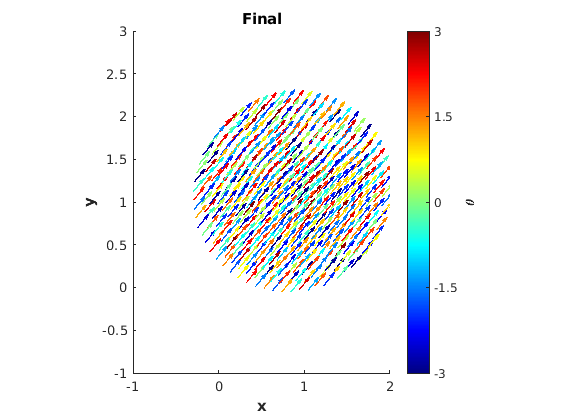}
\caption{Static async state \\ $(J,K) = (0.1,-1)$}
\end{subfigure}
\begin{subfigure}{0.32\textwidth}
\centering
\includegraphics[width=\linewidth]{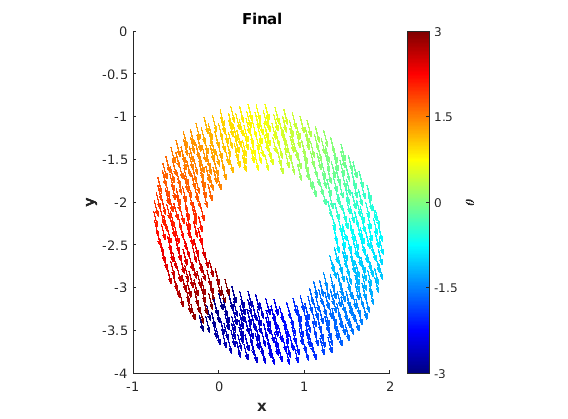}
\caption{Static phase wave state \\ $(J,K) = (1,0)$}
\end{subfigure}
\caption{Scatter plots of the three stationary states in the $(x,y)$ plane with $N = 500$ swarmalators for $T = 300$ time units.}
\end{figure}
Next, the moving states are also present.
\begin{figure}[H]
\centering
\begin{subfigure}{0.32\textwidth}
\centering
\includegraphics[width=\linewidth]{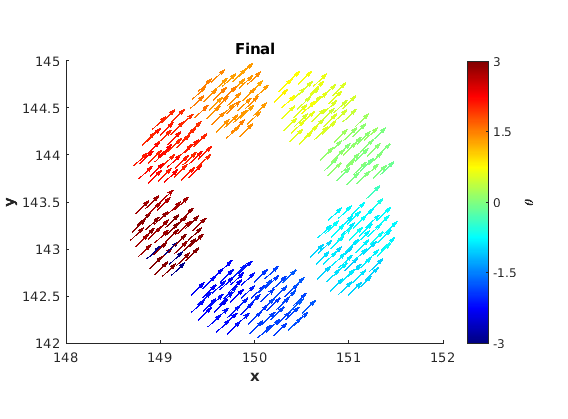}
\caption{Splintered phase wave state \\ $(J,K) = (1,-0.1)$}
\end{subfigure}\hfill
\begin{subfigure}{0.32\textwidth}
\centering
\includegraphics[width=\linewidth]{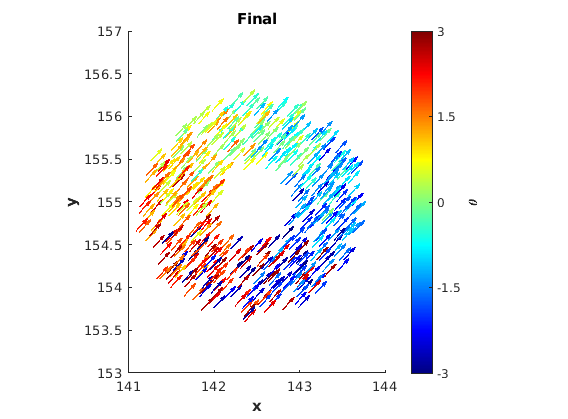}
\caption{Active phase wave state \\ $(J,K) = (1,-0.75)$}
\end{subfigure}
\caption{Scatter plots of the two active states in the $(x,y)$ plane with $N = 500$ swarmalators for $T = 300$ time units.}
\end{figure}
\section{Conclusion}\label{sec:conclude}
We have examined the dynamics of two types of systems which give concurrent emergence of clustering, flocking and synchronization.  The two systems describe orientable swarmalators with a local or global neighborhood.  These agents have spatial, phase, and orientation degrees of freedom which grants them the ability to flock (common velocity) and cluster (spatial pattern), as well as synchronize (common phase).  The position/velocity  (or heading)/phase are coupled through a set of interaction functions.  These models would be the first step in studying the behavior of real-world swarming oscillators (for example synchronized firefights).  We also look to biological systems for motivation to choose the radius of interaction, based on sensory limitations.  There are also a number of theoretical questions left to open, to be answered in future work.	
\bibliographystyle{abbrv}
\bibliography{sync}			
\end{document}